
\input eplain
\input amstex
\input epsf.tex
\documentstyle{amsppt}

    \pageheight{42pc}
    \pagewidth{27pc}


\chardef\oldatcatcode=\the\catcode`\@
\catcode `\@=11
\def\logo@{}
\catcode `@=\oldatcatcode


\NoBlackBoxes

\def\Z{\Bbb Z}

\def\C{\Bbb C}
\def\Q{\Bbb Q}
\def\Qbar{\overline{\Bbb Q}}
\def\R{\Bbb R}

\def\M{\Cal M}

\def\K{\Cal K}
\def\L{\Cal L}

\topmatter
\title{
Identities between Mahler measures
}\endtitle

\author Fernando Rodriguez-Villegas \endauthor
\affil{Department of Mathematics \\
  University of Texas at Austin \\
 Austin, TX 78712}\endaffil 
\address{Department of Mathematics, 
University of Texas at Austin}\endaddress 
\address Preprint\endaddress
\date{
September 2000
}\enddate

\thanks{
Support for this work was provided in part by a grant of the NSF 
and by a Sloan Research Fellowship.}
\endthanks

\email
villegas\@math.utexas.edu
\endemail

\endtopmatter

\document

\magnification=\magstep1

\baselineskip=14pt

The purpose of this short note is to give a  proof of the
following identity between (logarithmic) Mahler measures
$$
m(y^2+2xy+y-x^3-2x^2-x)=\frac57\,m(y^2+4xy+y-x^3+x^2)\,,
\tag{1}
$$
which is one of many examples that arise from the  comparison
of  Mahler measures and special values of $L$-functions
[Bo], [De], [RV]. Let us recall that the {\it logarithmic
Mahler measure} of a Laurent polynomial  $P \in \C 
[x_1^{\pm1}, \ldots , \;  x_n^{\pm1}]$
is defined  as
$$
m(P) = \int_0^1 \cdots \int_0^1 \, \log \left| P(e^{2\pi i \theta_1}
\, , \ldots , \; 
      e^{2\pi i\theta n}) \right| d\theta_1 \cdots d\theta _n \,.
\tag{2}
$$

The conjecture of Bloch--Beilinson [Be], [BG] for
elliptic curves predicts that both sides of (1) are rationally related
to $L'(E,0)$ (and hence to each other), where $E$ is the elliptic
curve of conductor $37$
$$
E:\quad y^2+y=x^3-x\,,
\tag{3}
$$
and $L(E,s)$ is its $L$-function. More precisely, we expect that the
 two numbers $a$ and $b$ defined by 
$$
\aligned
m(y^2+2xy+y-x^3-2x^2-x)&= a\, L'(E,0),\\ 
m(y^2+4xy+y-x^3+x^2) &= b\, L'(E,0)
\endaligned
$$
are rational. A proof of this fact is not without reach but will not
be attempted here, we will prove instead that $a/b=5/7$.

\vskip .5cm
\noindent
{\bf 1. Computing in $K_2(E)$}
\vskip .5cm
We first recall the definition of the group $K_2(A)$ of an elliptic
curve $A$. Given a field $F$ the group $K_2(F)$ can be defined as
$F^*\otimes F^*$ modulo the Steinberg relations $x\otimes (1-x)$ for
$x\neq 0,1$ in $F$.

Given a discrete valuation $v$ on $F$ with maximal ideal 
$\M$ and residue field $k$ we have the  {\it tame symbol} at $v$
defined by
 $$
 (x,y)_v\equiv (-1)^{v(x)v(y)}\frac{x^{v(y)}}{y^{v(x)}}
\mod \M,
 $$
which determines a homomorphism
$$
\lambda_v: K_2(F) \longrightarrow k^*
$$

For an elliptic curve $A$ defined over $\Q$  we let $K_2(A)$ be the
elements of $K_2(\Q(A))$ anihilated by all $\lambda_v$ with 
$v$ the valuations associated to $\Qbar$ points of $A$.

Our $E$ appears as a fiber in several of Boyd's families of
elliptic curves (see [Bo], [RV] for a discussion of these
families). For example, in its original form $y^2+y=x^3-x$, but  
also  as the two Weierstrass equations
$$
E_1: \quad y_1^2+4x_1y_1+y_1=x_1^3-x_1^2
\tag{4}
$$
and
$$
E_2:\quad y_2^2+2x_2y_2+y_2=x_2^3+2x_2^2+x_2\,.
\tag{5}
$$
It is easy to check that
$$
\aligned
x_1&=x-1\\
y_1&=y-2x+2
\endaligned
\tag{6}
$$
and
$$
\aligned
x_2&=x-1\\
y_2&=-x+y+1
\endaligned
\tag{7}
$$
give isomorphisms 
$$
E\simeq E_1, \qquad E\simeq E_2\,.
$$
It follows from [RV] therefore, that some integer multiple of  each of 
$$
\xi=\{x,y\}, \quad \xi_1=\{x_1,y_1\}, \quad \xi_2=\{x_2,y_2\}
$$
is in $K_2(E)$. 

The divisors of the six functions $x,y,x_1,y_1,x_2,y_2$ are supported
on $E(\Q)$, which is generated by the point $P$ with $x=0,y=0$. More
precisely, we have
$$
\aligned
(x)&= [P]+[-P]-2[O]\\
(y)&=[P]+[2P]+[-3P]-3[O]\\
\\
(x_1)&=[2P]+[-2P]-2[O]\\
(y_1)&=2[2P]+[-4P]-3[O]\\
\\
(x_2)&=[-2P]+[2P]-2[O]\\
(y_2)&=[2P]+2[-P]-3[O]
\endaligned
\tag{8}
$$
where $[O]$ denotes the point at infinity on $E$.

Given a pair of functions $f$ and $g$ on $E$ with divisors supported
on $E(\Q)$
$$
(f)=\sum_{n\in\Z} a_n[nP],
\qquad
(g)=\sum_{n\in\Z} b_n[nP]
$$
we define
$$
(f)\diamond(g)=\sum_{m,n}a_nb_m[(n-m)P]\,,
\tag{9}
$$
which we will view as an element of 
$$
\Z[E(\Q)]^-=\Z[E(\Q)]/\sim\,,
$$
where $\sim$ is the equivalence relation determined by
$$
[-nP]\sim -[nP], \qquad n\in \Z\,.
$$

We may and will represent elements of $\Z[E(\Q)]^-$ as vectors
$[a_1,a_2,\ldots]$ with $a_i\in \Z$ almost all zero where
$$
[a_1,a_2,\ldots] \qquad \longleftrightarrow \qquad \sum_{n=1}^\infty
a_n[nP] 
$$
In fact, we will only consider elements where $a_n=0$ for $n>6$ and
hence simply write $[a_1,\ldots,a_6]$.

We now compute
$$
\aligned
(x)\diamond (y)&= [1,2,-3,1,0,0]\\
(x_1)\diamond(y_1)&=[0,5,0,-4,0,1]\\
(x_2)\diamond(y_2)&=[-6,2,2,-1,0,0]\,.
\endaligned
\tag{10}
$$
On the other hand, we also find
$$
\aligned
(-y)\diamond(1+y)&=[-8,-7,8,1,0,-1]\\
(x-y)\diamond(1-x+y)&=[-9,5,-5,5,0,-1]
\endaligned
\tag{11}
$$
and verify easily that
$$
\aligned
7(x)\diamond(y)+(x_1)\diamond(y_1)&=-2(-y)\diamond(1+y)+(x-y)\diamond(1-x+y)\\
5(x)\diamond(y)+(x_2)\diamond(y_2)&=-(-y)\diamond(1+y)+(x-y)\diamond(1-x+y)\,.
\endaligned
\tag{12}
$$

\vskip .5cm
\noindent
{\bf 2. The regulator}
\vskip .5cm

Let
$$
r:\quad K_2(E) \longrightarrow \R
\tag{13}
$$
be the regulator map. It can be defined as follows. If $f,g$ are two
non-constant functions on $E$ with $\{f,g\} \in K_2(E)$ then
$$
r(\{f,g\})=\int_\gamma \eta(f,g),
\tag{14}
$$
where
$$
\eta(f,g)=\log|f|\; d\arg g - \log |g|\; d\arg f
\tag{15}
$$
and $\gamma$ is a closed path not going through poles or zeroes of $f$
or $g$ which generates the subgroup $H_1(E,\Z)^-$of
$H_1(E,\Z)$ where complex conjugation acts by $-1$, properly
oriented. The fact that the integral only depends on the homology
class of $\gamma$ is a consequence of $\{f,g\} \in K_2(E)$, see [RV]
for details. (However, note that in [RV] we inaccurately said $\gamma$
should generate the cycles fixed by complex conjugation; we take the
opportunity to correct this.) 

The regulator may also be expressed in terms of
the elliptic dilogarithm [BG], [Za] 
$$
\L: \quad E(\C) \longrightarrow \R\,.
$$
In our context, this works as follows. We extend it by linearity to
$\Z[E(\Q)]$ and since $\L$ is odd it actually gives a map
$$
\L: \quad \Z[E(\Q)]^- \longrightarrow \R\,.
\tag{16}
$$
 If $f,g$ are two non-constant functions on $E$ with divisors
supported on $E(\Q)$ and such that $\{f,g\} \in K_2(E)$ then
$$
r(\{f,g\})=c\;\L\left((f)\diamond(g)\right),
\tag{17}
$$
for some explicit non-zero constant $c$, which is not relevant
for our purposes. In particular, in the case that $g=1-f$ 
$$
\L((f)\diamond(1-f))=0\,.
\tag{18}
$$
The above discussion extends naturally to $K_2(E)\otimes \Q$, which
contains $\xi,\xi_1$ and $\xi_2$.

It follows from (12) therefore, that 
$$
\aligned
r(\xi_1)&=-7r(\xi)\\
r(\xi_2)&=-5r(\xi)\,.
\endaligned
\tag{19}
$$

\vskip .5cm
\noindent
{\bf 3. The regulator and Mahler's measure}
\vskip .5cm

In [RV] we showed that if $P_k(x,y)=0$ is one of Boyd's families of
elliptic curves and $k$ is such that $P_k$ does not vanish on the
torus $|x|=|y|=1$ then
$$
r(\{x,y\})=c_k\pi\;m(P_k)
\tag{20}
$$
for some nonzero integer $c_k$. We will now make this precise for
$$
P_k(x,y)= y^2-kxy+y-x^3+x^2\,.
$$
We consider the region $\K$ of $k\in \C$ such that $P_k$ vanishes
somewhere on the torus. It is the image of the torus under the
rational map 
$$
R: \qquad (x,y)\mapsto \frac{y^2+y-x^3+x^2}{xy}\,.
\tag{21}
$$
We can get a pretty good idea of what $\K$ looks like by graphing  the
image of a grid under $(\theta_1,\theta_2) \mapsto R(e^{2\pi i
\theta_1},e^{2\pi i \theta_2} )$. Dividing the square $0\leq
\theta_1<1, 0\leq \theta_2<1$ in $40$ equal parts we obtain
\midinsert
\botcaption{Figure 1} Region $\K$
\endcaption
\endinsert
It is not hard to verify directly that the boundary of $\K$ meets the
real axis at $k=-4$ and $k=2$.

If $k\notin \K$ then as $x$ moves counterclockwise on the circle
$|x|=1$ one root $y_1(x)$ of $P_k(x,y)=0$ satisfies $|y_1(x)|<1$ and
the other $y_2(x)$ satisfies $|y_2(x)|>1$ and in particular $y_1(x)$
and $y_2(x)$ do not meet. To see this, note that when $x=1$ the roots
are $0$ and $k-1$. Hence, for $|k|$ large these roots are one inside
and the other outside the unit circle. The claim follows since the
roots depends continuously on $k$. We let $\sigma_k$ be the resulting
smooth closed path $(x,y_1(x))$ on the elliptic curve $E_k$ determined
by $P_k(x,y)=0$.

Using Jensen's formula we find that 
$$
m(P_k)= \frac{1}{2\pi i} \int_{\sigma_k}\log|y| \;\frac{dx}x
$$
and note that since $|x|=1$ on $\sigma_k$ we can  write this identity as
$$
m(P_k)=\frac1{2\pi}\int_{\sigma_k} \eta(x,y)\,.
\tag{22}
$$

 We now show that for real and $k\notin \K$ the homology class of
$\sigma_k$ generates $H_1(E_k,\Z)^-$. We complete the square and write
$P_k=(2y-kx+1)^2-f(x)$, where $f(x)=4x^3+(k^2-4)x^2-2kx+1$. 
The discriminant $\Delta(k)=k^4-k^3-8k^2+36k-1$ of $f$ has two real
roots $\alpha=-3.7996\ldots$ and $\beta=.3305\ldots$. Hence, for
$k<\alpha$ or $k>\beta$, $\Delta(k)>0$ and $f$ has three real roots
$e_1<e_2<e_3$. As $|k|$ increases the roots of $f$ tend to
$e_1=-\infty$ and  $e_2=e_3=0$ and by continuity the circle $|x|=1$
encircles $e_2$ and $e_3$ once. Since $f$ is
negative in the interval $e_2<x<e_3$ the period
$$
\int_{\sigma_k}\frac{dx}{2y-kx+1}
$$
is purely imaginary and our claim follows.

Combined with (14) and (22) this  proves that in fact 
$$
r(\{x,y\})=\pm 2\pi\, m(P_k), \qquad k \in \R, \quad k\notin \K\,.
\tag{23}
$$
By continuity (23) also holds for $k=-4$ and $k=2$, which are on the
boundary of $\K$. In particular, in the notation of \S 2, we obtain
the identity
$$
r(\xi_1)=\pm 2\pi\, m(y^2+4xy+y-x^3+x^2)\,.
\tag{24}
$$
A completely analogous analysis yields
$$
r(\xi_2)=\pm 2\pi\, m(y^2+2xy+y-x^3-2x^2-x)
\tag{25}
$$
(and again $k=-2$ is on the boundary of the corresponding set $\K$). 
Putting together (19), (24) and (25) (and a simple check for the
right sign) we obtain (1).

\vskip .5cm
{\bf Remarks} \ 1. \ We should point out that we do not expect
$m(y^2+y-x^3+x)$  to be rationally related to either side of (1) (and
numerically it indeed does not appear to be). The reason is that
$y^2+y-x^3+x$ vanishes on the torus and in fact $k=0$ is in the
interior of the region $\K$ corresponding to the Boyd family
$y^2-kxy+y-x^3+x$. Hence the analogue of (22) gives the integral of
$\eta(x,y)$ on a non-closed cycle.

2. \ One can prove in a similar way an identity relating either side
 of (1) with $m(y^2+2xy+y-x^3+x^2)$.

\enddocument